\documentclass{amsart}

\usepackage[T1]{fontenc}
\usepackage[utf8]{inputenc}
\usepackage{lmodern}
\usepackage{url}
\usepackage{hyperref}

\title[Explanations, Prompts, and Formalizations]{Explanations, Prompts, and Formalizations:
Arguments for New Norms in LLM-Enabled Mathematical Research}
\author{Axel Boldt}
\address{Department of Mathematics and Statistics,
Metropolitan State University, Saint Paul, Minnesota, U.S.A.}
\email{Axel.Boldt@metrostate.edu}
\date{6 August 2026}

\begin{document}

\maketitle

\begin{abstract}
As several mathematical conjectures have recently been settled using large language models (LLMs), the mathematical community has formulated norms and recommendations regarding the publishing of such results. These norms do not cover the disclosure of the prompts and precise software setup used to obtain those results, nor do they require that results be formalized in a manner that allows for machine verification. I argue that both of these are essential. In addition, since LLM-obtained results may be hard to understand, human authors have the responsibility to invent intuitive explanations.
\end{abstract}

\section{Recent LLM-supported developments in mathematics}

The mathematical community has for several years worked to translate known mathematical definitions, theorems, and proofs into formal computer languages, thereby allowing proofs to be verified automatically and facilitating machine-assisted construction of new proofs. Perhaps the most advanced effort in this area is \texttt{mathlib}, a large library of definitions, theorems, and proofs in the \texttt{Lean} language.\cite{ullrich2024} While formalizing a given mathematical result is a laborious and nontrivial task, recently some large language models (LLMs) have made advances in this realm.\cite{hariharan2026}

Beginning in late 2025, LLMs became significantly more capable of proving or disproving mathematical conjectures outright. Such a system, when given an appropriate prompt, can search the internet for prior work and proceed to construct, albeit sometimes incorrectly, proofs or disproofs written in the mathematical vernacular. Several minor conjectures were settled in this manner, many of them formulated by the prolific mathematician Paul Erd\H{o}s (1913--1996).

In order to get a better understanding of the mathematical capabilities of LLMs, in February 2026 the First Proof project released a batch of ten research-level mathematical questions that had previously been solved by experts, but whose solutions had not yet been published.\cite{abouzaid-first-proof} LLM solutions were collected, and many of them were found to be unsatisfying or defective. A second batch of similar problems was released in June; this time, it was guaranteed that the LLMs worked entirely autonomously, with no human help, and each solution was evaluated by several independent blinded experts. The prompts and harnesses---a harness is software that interacts with and guides an LLM, for instance by providing it with further prompts, data, or tools based on its prior output---were fully published. Four systems were tested, each on the ten problems; seven of the ten problems received at least one acceptable solution.\cite{abouzaid-second-batch}

A breakthrough occurred in May 2026, when an LLM settled a major problem in combinatorial geometry, Erd\H{o}s's unit-distance conjecture. This conjecture was widely believed to be true; the LLM disproved it by constructing a counterexample.\cite{openai-unit-distance} The conjecture concerns \(n\) points in the plane and claims a bound on the maximal number of pairs of points that have mutual distance \(1\). The counterexample used an innovative construction from algebraic number theory. In an accompanying article, a number of experts explained the proof and the result's significance.\cite{alon-remarks} The LLM's ``Chain-of-Thought'' document was also released; this is akin to scratch paper used by the LLM for brainstorming while constructing the proof. Neither the exact version of the LLM nor the exact prompt nor any harness, if one was used, was revealed. A formalization of the result in \texttt{Lean} was constructed with LLM assistance one month later.\cite{lean-erdos-unit}

The cycle double cover conjecture, a famous and longstanding open problem in graph theory, was proved true by an LLM in July 2026. The short proof built on prior work and employed one new and essential trick. It was released along with the exact prompt used to create the proof, as well as a formalization in \texttt{Lean} of the prior work and of the proof.\cite{openai-cdc}

Later in the same month, the important Jacobian conjecture in algebraic geometry was disproved by an LLM. The short counterexample was published on social media, without any elaboration.\cite{alpoge-jacobian} Shortly thereafter, other authors attempted to explain the counterexample, with LLM assistance.\cite{tao-jacobian,lou-jacobian}

\section{Reactions and evolving norms}

Since solving problems, proving theorems, and constructing counterexamples have traditionally been among the main activities of research mathematicians, it comes as no surprise that the mathematical community had to formulate a response. In a widely discussed blog post, David Bessis criticizes Geoff Hinton's view that mathematics, like chess and Go, is a closed system with rules, ready for the taking by AI.\cite{bessis-theorem-economy} He approvingly cites prominent mathematicians Bill Thurston and Terence Tao, who argue that the true task of mathematicians is to create human understanding, not just produce definitions, theorems, and proofs.

The widely endorsed Leiden Declaration of June 2026 records a number of norms and expectations regarding the use of LLMs by authors and publishers of mathematical works.\cite{leiden-declaration} For our purposes, the most significant of these are: (a) LLM use should always be disclosed; (b) authors should make an effort to track down references to prior work even if the LLM did not cite those works; and (c) formalizations should be provided if feasible and appropriate.

Regarding LLM-produced formalizations of new proofs, it has been pointed out that these, while desirable, are often not provided in a form that allows incorporation into \texttt{mathlib}, which follows conventions and norms intended to facilitate reuse of the results in further formalizations. The process of turning raw formalizations into useful additions to \texttt{mathlib} is known as canonization; at the moment, this task is almost exclusively performed by humans.\cite{kontorovich-canonization}

\section{Demand on authors: explanations}

If we endorse the common view that the foremost task of a mathematician is to facilitate and enlarge mathematical understanding, it becomes obvious that mathematical results obtained by LLMs should be accompanied by human-written explanations. LLMs often do not properly cite all relevant prior work, even work that may have been part of their training sets and may have contributed to the final output. Furthermore, solutions written by LLMs will often explain routine steps of a derivation in exhausting detail while glossing over the new and interesting parts, as was pointed out by several referees of the First Proof second batch.\cite{abouzaid-second-batch} This is likely a consequence of the LLM's training: the routine steps, unlike the innovative ones, will have occurred repeatedly in its training set. It is clear that the human mathematician must correct this disproportion in a pedagogical manner. Lengthy calculations can sometimes be summarized or explained in a conceptual manner or illustrated by suggestive figures; it is the task of the human to provide these if they are not already present.

LLM-provided proofs commonly contain subtle mistakes, some small and some fatal. This may be a consequence of their training procedure: they are designed to produce text that statistically matches proofs---text that looks like a proof---but such text need not {\em be} a proof. As the Leiden Declaration also points out, the human author must accept full responsibility for the correctness of the published results,\cite{leiden-declaration} which necessitates a complete understanding of all arguments.

At the current stage of LLM development, human experts retain a better overview of a field, can explain the significance of a result, point out connections to other problems and developments, and suggest future research directions. All of these ought to be required of human reporters of LLM-derived results.

\section{Demand on authors: formalizations}

As mentioned above, LLM-supplied proofs can be, and often are, wrong, typically in very subtle ways. LLMs can also aid in the formalization of definitions, theorems, and proofs, so that the work may be checked automatically. I believe it is the responsibility of a mathematician reporting on LLM-obtained work to provide such a formalization.

To be sure, a formalization does not provide absolute certainty regarding the correctness of the work. There can be subtle mistakes in the formalization of the definitions and theorems, so that the formalized work does not precisely match the informal presentation. Furthermore, LLMs have been known to detect bugs in the \texttt{Lean} kernel and other automated proof checkers, allowing them to cheat and construct defective proofs that the proof checkers accept as correct.\cite{demoura-postmortem} Nevertheless, a formalized proof has a much higher likelihood of being correct and is easier to check than a complicated, lengthy, and novel informal argument.

As a best practice, it would be desirable that the formalization be canonized and follow the conventions of \texttt{mathlib} so that it can easily be incorporated into that library, accelerating the growth of that increasingly important resource. However, producing such a canonized formalization requires quite specialized knowledge, so I would not demand it of authors.

\section{Demand on authors: prompts and harnesses}

Several of the reports about LLM-produced mathematical work that I described in Section~1 lack crucial information. The precise version of the LLM is often not mentioned; sometimes it is described only as an ``internal model.'' The precise prompt or prompts given to the LLM  are often not given. Nor, if a harness was used, is the precise version of that harness mentioned. I consider this deplorable.

There is a danger that the ``art of prompting'' could turn into a secret art, known only to insiders. This would be reminiscent of sixteenth-century mathematicians hiding their solution to the cubic equation, or of alchemists keeping their experiments secret. In order to maximally expand human understanding, we should not only explain the mathematics but also how the mathematics was discovered. This transparency will allow others to learn.

A laudable example is the proof of the cycle double cover conjecture: the entire two-page prompt was published.\cite{openai-cdc} It was undoubtedly written by experts who knew about the pitfalls of the conjecture, and it also contained several new and interesting techniques, such as having a number of independent agents pursue different approaches in parallel and having adversarial agents double-check one another's work. This is crucially important information that can help many other mathematicians in their work. The version of the LLM, however, was not revealed.

It is true that, even if the LLM version, prompt(s), and harness are specified, the output of a system is not fully deterministic and the experiment is therefore not completely reproducible. But given that information, the experiment can be repeated several times and the outputs compared. The situation is akin to that in chemistry, where researchers are expected to describe their experimental setup in as much detail as possible, knowing that an element of chance will always remain.

\section{Objections and rationale}

Authors of articles reporting on LLM-produced work might object that requirements such as those outlined in this article, which are not normally imposed on other mathematicians, are unfair. I would argue that these authors already do substantially less work than other mathematicians and should face a higher ``burden of proof.'' This is, in part, a defense of the refereeing and publishing pipeline against a flood of quickly produced, automatically generated work of low quality. It also follows from the overarching goal of maximizing all aspects of human understanding.

\end{document}